\documentclass{article}
\usepackage{amsthm}
\usepackage{amsfonts}
\usepackage{amssymb}
\usepackage{amsmath}
\usepackage{makeidx}
\usepackage[english]{babel}
\usepackage{xypic}
\newtheorem{teo}{Theorem}[section]
\theoremstyle{definition}
\newtheorem{defi}[teo]{Definition}
\theoremstyle{remark}
\newtheorem{rem}[teo]{Remark}
\theoremstyle{remark}

\theoremstyle{theorem}
\newtheorem{lemma}[teo]{Lemma}
\theoremstyle{theorem}
\newtheorem{prop}[teo]{Proposition}
\newtheorem{cor}[teo]{Corollary}

\newcommand{\lgra}{\longrightarrow}

\newcommand{\rd} {\mathrm{rdim\ }}
\renewcommand{\deg}{\mathrm{deg }}
\renewcommand{\dim}{\mathrm{dim}}

\newcommand{\p}{\mathfrak{p}}
\newcommand{\m}{\mathfrak{m}}
\newcommand{\bb}{\mathfrak{b}}

\newcommand{\M}{\mathfrak{M}}
\newcommand{\aid}{\mathfrak{a}}
\newcommand{\annu}[2]{(0:_{#1}{#2})}
\newcommand{\gr}[2]{\mathrm{gr}_{#1}{#2}}
\newcommand{\biadeg}[2]{\underline{\mathrm{arithdeg}}_{#1}(#2)}
\newcommand{\locadeg}[3]{\underline{\mathrm{arithdeg}}_{#1}(#2,#3)}
\newcommand{\adeg}[2]{\mathrm{arithdeg}_{#1}(#2)}

\newcommand{\gmult}[3]{\underline{\mathrm{e}}_{#1}({#2},{#3})}
\newcommand{\ee}[2]{\underline{\mathrm{e}}_{#1}({#2})}
\newcommand{\e}[2]{\mathrm{e}_{#1}({#2})}
\newcommand{\emult}[3]{\mathrm{e}_{#1}({#2},{#3})}
\newcommand{\ext}[4]{\mathrm{Ext}_{#1} ^{#2}({#3},{#4})}

\newcommand{\ass}[1]{\mathrm{Ass}({#1})}
\newcommand{\Ho}[2]{\mathrm{H}^{0}_{#1}({#2}_{#1})}
\newcommand{\ann}[1]{\mathrm{Ann}{(#1)}}
\newcommand{\fg}{finitely generated }

\begin{document}
  \title{Arithmetic degree and associated graded modules}  
  \author{Natale Paolo Vinai \\ Dipartimento di Matematica, Universit\`a di Bologna,\\ Piazza di Porta San Donato 5, \\40126 Bologna, Italy\\ \texttt{vinai@dm.unibo.it} \\    \date{}   } 
\maketitle



 \begin{abstract} 
We prove that the arithmetic degree of a graded or local ring $A$ is bounded above by the arithmetic degree of any of its associated graded rings with respect to ideals $I$ in $A$. In particular, if Spec$(A)$ is equidimensional and has an embedded component (i.e., $A$ has an embedded associated prime ideal), then the normal cone of Spec$(A)$ along V$(I)$ has an embedded component too. This extends a result of W.~M.~Ruppert about embedded components of the tangent cone.
 \end{abstract}
\section{Introduction}

It is well-known that many homological properties of a Noetherian commutative ring $A$ can be recovered from the properties of its associated graded ring $\gr{I}{A} = \displaystyle{\oplus_{n\geq 0}} I^{n}/I^{n+1} $ with respect to some ideal $I$ in $A$. For example, if $(A,\m)$ is a local ring and $\gr{I}{A}$ is regular (complete intersection, Gorenstein, Cohen-Macaulay) then $A$ is so, and $\e{}{A} \leq \e{}{\gr{I}{A}}$ where $\mathrm{e}$ denotes the Samuel multiplicity with respect to the unique maximal ideal of $A$ or the unique homogeneous maximal ideal of $\gr{I}{A}$, respectively.
The present paper is motivated by a more recent result of W. M. Ruppert \cite{Ru} who proved that if the ring $\gr{\m}{A}$ has no embedded associated prime ideals, then $A$ has also no embedded associated prime ideals. We extend and refine this result by showing that in fact the arithmetic degree of $A$ is bounded above by the arithmetic degree of $\gr{I}{A}$, where $I$ is an arbitrary ideal of the local ring $(A,\m)$ (see corollary of Theorem \ref{Ineq}):
\begin{center}
$\adeg{}{A} \leq \adeg{}{\gr{I}{A}}$.
\end{center} 
The idea of the  arithmetic degree goes back to Hartshorne \cite{Har} and the name ``arithmetic degree'' has been introduced by D. Bayer and D. Mumford in \cite{BM}. The arithmetic degree is a measure of the complexity of a homogeneous ideal in a polynomial ring, which refines the classical degree.
Whereas the degree takes into account only the top-dimensional primary components of the ideal, the arithmetic degree involves also lower dimensional (both isolated and embedded) components.

Our results and their proofs require an extension of the notion of
arithmetic degree to modules over a local ring (or a standard graded
algebra) given by W. V. Vasconcelos \cite{Vas}, p.~223, and the study
of the arithmetic degree of bigraded modules, see section 2. 

\begin{section}{Bigraded rings and modules}

Let $A$ be a finitely generated bigraded algebra over an Artinian local ring $A_{(0,0)}$.
If the generators of $A$ are of bidegree $(1,0)$ and $(0,1)$ we will say that $A$ is a standard bigraded algebra. We will mostly work in this case.

Indeed, in this setting given a finitely generated $A$-module $M$, we can consider its bigraded Hilbert function which turns out to be a polynomial (see e.g. \cite{ Trung1}) in two variables.
As in the graded case some important informations are encoded in this polynomial, e.g. (with some modifications) the dimension of the module. 
Let us recall some properties of the dimension.

If $x_{1}, \dots , x_{n}$ are the generators of $A$ over $A_{(0,0)}$ of bidegree $(1,0)$ and $y_{1}, \dots , y_{m}$ those of bidegree $(0,1)$, then we denote by $A_{+}$ the ideal $(x_{1},\dots,x_{n})\cap (y_{1}, \dots, y_{m})$.
An ideal which does not contain $A_{+}$ will be called {\em{relevant}}.
If $M$ is an $A$-module, $N$ a submodule of $M$ and $I$ an ideal in $A$, let
 \begin{center}$ (N:_{M}I^{\infty})= \{ x \in M \mid$ there is a positive integer n such that $ x I^{n} \subseteq N \}.$ \end{center}
If $M$ is $A$ we write simply $(N:I^{\infty})$.

 \begin{defi}

The relevant dimension of $A$ is

\begin{center}

$\rd A = \dim A/(0:A_{+}^{\infty}).$

\end{center}

For a bigraded finitely generated $A$-module $M$, we set

\begin{center}

$\rd M = \rd A/\ann{M}.$

\end{center}

\end{defi}

\begin{rem}

We will need the more useful charaterization

\begin{center}
$\rd M = \mathrm{max} \{ \dim A/\p \mid \p $ relevant bihomogeneous prime ideal containing $ \ann{M} \}$.
\end{center}

\end{rem}

As the notion of dimension, relevant dimension behaves well with respect to short exact sequences. With the same proof as in the graded case, we obtain the following result:

\begin{lemma} \label{maxrd}
Let 
\begin{center}
$0 \lgra M' \lgra M \lgra M'' \lgra 0$
\end{center}
be a short exact sequence of finitely generated bigraded $A$-modules.
Then
\begin{center}
$\rd M = \mathrm{max}\{ \rd M', \rd M'' \}.$
\end{center}
\end{lemma}

We have given a notion of relevant prime ideal consistent with the notion of relevant dimension. We need to extend this idea to modules with a new definition, which is still close to the graded one.

\begin{defi}
A \fg  $A$-module $M$ is {\em{relevant}} if there exists an irredundant primary decomposition
\begin{center}
$0_{M} = \cap M_{i}$ ;   $\ass{M/M_{i}} = \{\p_{i}\}$
\end{center}
such that  $\p_{i}$ is relevant for at least one $i$.
\end{defi}

\begin{rem} \label{RelProp}
Relevant modules $M$ have the following properties which can be easily shown:
\begin{itemize}
\item[(i)] $\mathrm{Ann}(M) \nsupseteq A_{+}$;
\item[(ii)] If $N \supset M$ then $N$ is relevant;
\item[(iii)] There exists a filtration of $M$
 \begin{center}
$ 0 = M _{0} \subset M_{1} \subset \dots \subset M_{n} = M$ 
\end{center}
such that $M_{i}/M_{i-1} \simeq A/\p_{i}$ with $\p_{i}$ bihomogeneous prime ideal and $\p_{1}$ is a relevant ideal;

\item[(iv)] In the above filtration all the $M_{i}$'s are relevant.
\end{itemize}
\end{rem}

\noindent Another property of a relevant module is the following generalization of a result of N.V.~Trung \cite{Trung1}.

\begin{prop} \label{Hilb}
Let $M\neq0$ be a relevant finitely generated  $A$-module and $P_{M}(r,s)$ its Hilbert polynomial. Then $\deg P_{M}(r,s) = \rd M -2$.
\end{prop}

\begin{proof}
The proof is by induction on the length $n$ of the filtration (iii) in Remark \ref{RelProp}.
Since $M \neq 0$ then $n > 0$ .
If $n = 1$ then $ M \simeq A/\p$, thus $M$ satisfies the hypothesis of Theorem 1.7 of \cite{Trung1} which is our thesis for ideals.

Suppose $n>1$ and that the assertion holds for $n-1$. From the filtration we have

\begin{center}

$0 \lgra M_{n-1} \lgra M \lgra A/\p \lgra 0$

\end{center}

\noindent 
and since the Hilbert polynomial is additive, $P_{M}(r,s) = P_{M_{n-1}}(r,s) + P_{A/\p}(r,s)$.
By the previous remark $M_{n-1}$ is relevant and so by induction its relevant dimension is the degree of its Hilbert polynomial plus two. Our proposition follows by induction on $n$ and Lemma \ref{maxrd}. 
\end{proof}

\begin{subsection}{Multiplicities of bigraded modules}

In this section we recall some facts on Hilbert functions and multiplicities of bigraded modules.

\begin{defi}
Let $H\colon \mathbb{Z}^{2} \lgra \mathbb{N}$ be a numerical function.We define inductively the differences $\Delta^{(m,n)}H(i,j)$:

\begin{itemize}
\item[] $\Delta^{(1,0)}H(i,j) = H(i,j) - H(i-1,j)$,
\item[] $\Delta^{(0,1)}H(i,j) = H(i,j) - H(i,j-1)$,
\item[] $\Delta^{(1,1)}H(i,j) = \Delta^{(1,0)}\Delta^{(0,1)}H(i,j)$,
\item[] $\Delta^{(m,n)}H(i,j) =  \left\{ \begin{array}{ll}
        \Delta ^{(m-n,0)} \Delta^{(n,n)}H(i,j) & \textrm{if } m\geq n \\
        \Delta ^{(0,n-m)} \Delta^{(m,m)}H(i,j) & \textrm{if } n\geq m .
        \end{array} \right. $
\end{itemize}
\end{defi}
\noindent
Composition of differences is well defined and 
\begin{center}
$\Delta^{(r,s)}\Delta^{(m,n)}H(i,j)=\Delta^{(r+m,s+n)}H(i,j).$
\end{center}
In particular this definition can be applied to numerical polynomials. It is well known that a numerical polynomial in two variables can be written in the form
\begin{center}
$P(m,n) = \displaystyle{\sum_{i+j \leq d;i,j \geq 0}}  a_{i,j} {m \choose i}{n \choose j}  $
\end{center}
with integers $a_{i,j}$ and by computation we obtain
\begin{center}$\Delta^{(t,s)}P(m,n)=\displaystyle{\sum} a_{i,j}{m-t \choose
 i-t}{n-s \choose j-s}$. \end{center}

Given a bigraded ring $A$ and a bigraded $A$-module $M$, the ``right'' polynomial to be considered to define multiplicity is the sum transform of the Hilbert polynomial as shown in \cite{AMa}.
We recall here some results of \cite{AMa}.
Let $h(i,j) = $ length $_{A_{(0,0)}} (M_{(i,j)})$  be the Hilbert function of $M$, and $h^{(1,0)}(i,j)=\sum_{u=0}^{i} h(u,j)$ the sum transform of $h$ with respect to the first variable.
Then we have still to consider the sum transform with respect to the second variable:
$h^{(1,1)}(i,j)=\sum_{v=0}^{j}h^{(1,0)}(i,v)$.
It is clear that, like for the Hilbert function, also this sum transform function becomes a polynomial for $i,j$ sufficiently large.
\noindent We observe that working with the double sum transform $h^{(1,1)}(i,j)$ has the same effect as adding two new variables $x,y$ to the bigraded ring, one of bidegree $(1,0)$ and one of bidegree $(0,1)$ (see Proposition 1.2 of \cite{AMa}) and working with $h(i,j)$. For an $A$-module $M$ this means that it has to be replaced by the $A[x,y]$-module $M[x,y] \simeq M \otimes_{A} A[x,y]$, i.e. 
$$ h_{M}^{(1,1)}(i,j) = h_{M[x,y]}(i,j)$$
The passage from $A$ to $A[x,y]$ and from $M$ to $M[x,y]$, by the point of view of the associated prime ideals, does not change anything since there is a one-to-one correspondence
\begin{center}
$\ass{M} \leftrightarrow \ass{M[x,y]}$
\end{center}
but the $A[x,y]$-module $M[x,y]$ is relevant. Furthermore, the relevant dimension agrees with the Krull dimension for a relevant module.
Since from now on we will always work with the sum transform of Hilbert functions, we can remove the hypothesis ``relevant'' from all our results. 

We want to define a multiplicity symbol $\ee{q}{M}$ which is additive on short exact sequences and that generalizes the graded multiplicity symbol.
Let $M$ be a bigraded \fg $A$-module of dimension $d$. We consider the double sum transform of its Hilbert polynomial $P_{M}^{(1,1)}(i,j)$ and the binomial expression $P_{M}^{(1,1)}(i,j) =  \displaystyle{\sum_{i+j\leq d}} c_{i,j} {m \choose i}{n \choose j}$.
\begin{defi} Under the above hypothesis, we set
$$\ee{q}{M}=\gmult{q}{A_{+}}{M} = \left\{\begin{array}{ll} (c_{0,d},\dots,c_{t,t-d},\dots,c_{d,0}) & q=\dim M \\
                          \underline{0}        & \textrm{otherwise}
                \end{array}\right. .$$
By $\underline{0}$ we mean the vector of $q+1$ components all equal to $0$. 
\end{defi}
\begin{rem}
If $M=A$ then this generalized multiplicity is the vector composed by the $c_{k}(A)$'s of \cite{AMa}, i.e. $(\ee{\rd{A}}{A})_{k} = c_{k}(A)$.
\end{rem}
The preceding definition is inspired by the classical definition of multiplicity in a local ring, which we recall here.
\begin{defi}
Let $(A,\m)$ be a Noetherian local ring, $M$ a \fg $A$-module. For an integer $i$ we define:
$$\e{i}{M} =  \left\{\begin{array}{ll} \e{}{M} & i=\dim M ;\\
                                   0        & \textrm{otherwise}.
                \end{array}\right.$$
\end{defi}
We can use differences in two variables to compute the generalized multiplicity symbol. In fact we obtain by induction on the dimension the following lemma.
\begin{lemma}\label{leadcoef}
Let $M \neq 0$ a bigraded \fg $A$-module and $(t,s) \in \mathbb{Z}^{2}$ such that $t+s = d = \dim M$. Then writing the sum transform of the Hilbert polynomial of $M$ in the form $$P_{M}^{(1,1)}(m,n) = \sum c_{i,j} {m \choose i}{n \choose j},  $$
we have
\begin{center}
$\Delta^{(t,s)}P_{M}^{(1,1)}(m,n) = c_{t,s}.$
\end{center}
\end{lemma}
As we required, these vectors are additive on short exact sequences.
\begin{prop}
Let
\begin{center}$0 \lgra M^{'} \lgra M \lgra M^{''} \lgra 0 $ \end{center}
be a short exact sequence of bigraded \fg $A$-modules such that $\dim
M \leq q $. \\
Then \begin{center} $\ee{q}{M}=\ee{q}{M^{'}} + \ee{q}{M^{''}}$.\end{center}

\end{prop}

\begin{proof}
We can assume $\dim M = q$.
By the additivity of the Hilbert polynomial on short exact sequences, $P^{(1,1)}_{M}(m,n)=P^{(1,1)}_{M^{'}}(m,n)+P^{(1,1)}_{M^{''}}(m,n)$. If we choose $s,t$ such that $s+t=q$, then we have
$$\Delta^{(t,s)}P^{(1,1)}_M(m,n)=c_{t,s}=\Delta^{(t,s)}P^{(1,1)}_{M^{'}}(m,n)+
\Delta^{(t,s)}P^{(1,1)}_{M^{'}}(m,n).$$ Taking into account that $\dim M = \mathrm{max}\{\dim M^{'}, \dim M^{''}\}$ and applying Lemma \ref{leadcoef}, the result follows.
\end{proof}

The standard application of the last proposition is:

\begin{prop}\label{add} {\rm{(see e.g. \cite{Eis})}}
Let $M$ be as above and $\ell(-)$ denote the length of a module on the appropriate ring. Then \begin{center} $\ee{q}{M} =$
$\displaystyle{\sum_{\dim A/\p = \dim M}} $ $\ell(M_{\p})\ee{q}{A/\p},$
\end{center}
where $\p$ runs over all highest-dimensional associated prime ideals of $M$.
\end{prop}

\end{subsection}

\begin{subsection}{Arithmetic Degree}
\noindent Now we have almost all the necessary notions to introduce the arithmetic degree for bigraded modules. 
We recall here Vasconcelos's definition of arithmetic degree \cite{Vas}, p.~223 for graded or local rings, since it will give us the basic idea to extend it to the bigraded setting.

\begin{defi}
Let $A$ a graded (or local) ring, $M$ a \fg $A$-module.  With $\ell(-)$ we denote the length of a module on the appropriate ring. For $i \geq 0$ we define:

\begin{center}
$\adeg{i}{M} = \displaystyle{\sum_{\dim A/\p = i }} \ell (\Ho{\p}{M}) \e{i}{A/p}. $
\end{center} 

\end{defi}
Let $A$ be a bigraded ring and $M$ an $A$-module.
To define the arithmetic degee for bigraded rings, we set
\begin{center}
$M_{\leq i} = ( x \in M \mid \dim (A/ \ann{x}) \leq i ) $
\end{center}
and $$M_{>i}=M/M_{\leq i}.$$

Actually these definitions are the same used in \cite{MVY} in the graded case. 

\noindent These two bigraded $A$-modules have the following  properties:
\begin{itemize}
\item[(i)] $\ass{M}_{\leq i} = \{ \p \in \ass{M} \mid \dim A/\p \leq i \}$; 
\item[(ii)] $\dim M_{\leq i} \leq i $;
\item[(iii)] $\ass{M}_{>i}= \ass{M} / \ass{M}_{\leq i}$;
\item[(iv)] If $M_{>i}\neq 0$ then $\dim M_{>i} > i$.
\end{itemize}

\begin{defi}
Let M be a \fg bigraded $A$-module. We define the $i$th bigraded arithmetic degree of $M$: 
\begin{center}
$\biadeg{i}{M} = \displaystyle { \sum_{\dim
(A/ \p) = i}} \ell((M_{\leq i})_{\p})\ee{i}{A/ \p}$
\end{center}
where $\p$ runs on the associated prime ideals of $M$.
\end{defi}

\noindent At first sight, this definition seems to be different from Vasconcelos's one \cite{Vas} p.~223, but in virtue of the next lemma the length of $(M_{\leq i})_{\p}$ is just the length multiplicity of $\p$ with respect to $M$ used in \cite{Vas}.

\begin{lemma} \label{LocCohom}
Let $M$ be a \fg $A$-module and $\p \in \ass{M}$ such
that $\dim A/ \p = i$. Then
\begin{center}
$H_{\p}^{0}(M_{\p}) \simeq (M_{\leq i})_{\p}$.
\end{center}
\end{lemma}

\begin{proof}
Since $(M_{\leq i})_{\p}$ is of finite length and $\Ho{\p}{M}$ is the largest submodule of finite length of $M_{\p}$ we have only to
prove that $H^{0}_{\p}(M_{\p}) \subseteq (M_{\leq i})_{\p}$. 
Take an element $x \in H^{0}_{\p}(M_{\p})$ such that
$x=\frac{m}{s}$ with $m \in H^{0}_{\p}(M)$, $s \in A - \p$. Then there exists $n \in \mathbb{N}$ such that $\p ^{n}m=0$ which means that $\sqrt{(0:_{A}
m)} \supseteq \p.\;\;$ It follows that \begin{center} $\dim {A/(0:_{A} m)} =
\dim {A/\sqrt{(0:_{A} m)}} \leq \dim {A/\p} = i.$ \end{center} Hence
$m \in M_{\leq i}$ and $x \in (M_{\leq i})_{\p}$
\end{proof}

\noindent Our definition has a historical motivation. It is quite closer to Hartshorne's one in \cite{Har}, where to notion of arithmetic degree was introduced for the first time. Indeed in \cite{Har} the arithmetic degree was defined in a form similar to that described by in the following lemma, which is a standard charaterization of the arithmetic degree (see e.g. \cite{STV} or \cite{MVY}).
If $H(m,n)$ is a numerical function, we denote by $\underline{\Delta}^{i}(H(m,n))$ the vector 
$$(\Delta^{(0,i)}H(m,n),\dots,\Delta^{(t,i-t)}H(m,n),\dots,\Delta^{(i,0)}H(m,n)).$$

\begin{lemma}
Let $M$ be a \fg bigraded $A$-module. Then:
\begin{center}$ \biadeg{i}{M} = \underline{\Delta^{i}}(P_{M}(m,n) -
              P_{M > i}(m,n))$   for all $ (m,n) \in \mathbb{Z}^{2} .$
\end{center}
\end{lemma}
\begin{proof}
Consider the short exact sequence
\begin{center}
$0 \rightarrow M_{\leq i} \rightarrow M \rightarrow M_{>i}
\rightarrow 0$.
\end{center}
By the linearity of the Hilbert function and the properties of
differences in two variables we get the result.
\end{proof}

\noindent Lemma \ref{LocCohom} has the important consequence \ref{Ext}, since it shows how to compute the arithmetic degree for standard bigraded rings using computer algebra programs (e.g. Reduce \cite{Red} using the package Segre \cite{AAl}, Macaulay \cite{Mac} and CoCoA \cite{Cocoa}) which compute Ext modules.

\begin{prop}\label{Ext}
In the above setting, the following holds:
\begin{center}
$\biadeg{i}{M} = \ee{i}{\ext{A}{\dim A - i}{M}{A}}.$
\end{center}
\end{prop}
\begin{proof}
It follows easily by  Lemma \ref{LocCohom} and local duality (see e.g. \cite{BroSh}).
\end{proof}

\end{subsection} 

\end{section}

\begin{section}{Generalized arithmetic multiplicity of a local ring}

 Now we shall define the arithmetic multiplicity relative to a proper ideal for a module on a local ring.
From now on $(A,\m)$ will be a Noetherian local ring.
\noindent For an $A$-module $M$ we will denote $\gr{\m}{\gr{I}{M}}$ by $GG(M)$, considered as a bigraded module over  the bigraded ring $GG(A)=\gr{\m}{\gr{I}{A}}$.
The bigraduation is given by
\begin{center}
$GG(M)_{(i,j)} = {\mathrm gr}_{\m}^{i} \mathrm{gr}_{I}^{j}{M}_{(i,j)}=(\m^{i}I^{j}+I^{j+1})M/ (\m^{i+1}I^{j}+I^{j+1})M $.
\end{center}

Extending to the case of modules a definition of \cite{AMa}, the multiplicity of an arbitrary ideal $I$ with respect to a finitely generated  $A$-module $M$ is defined to be the multiplicity of the bigraded module $GG(M)$.
\begin{defi}
Let $(A,\m)$ be a Noetherian local ring, $I$ an ideal in $A$ and $M$ a \fg $A$-module. For an integer $i$ we define
 $$\gmult{i}{I}{M} = \ee{i}{GG(M)}.$$
\end{defi}

In particular, if $I$ is $\m$-primary our definition agrees with the classical definition of multiplicity, as can be seen from the following proposition.

\begin{prop}{\rm{(see Corollary 2.4 of \cite{AMa})}} \label{clad}

\par If $I$ is $\m$-primary then 
\begin{center}
$\emult{ }{I}{A} = (\gmult{}{I}{A})_{0} = c_{0}(A)$
\end{center}
and 
\begin{center}
$(\gmult{}{I}{A})_{k}=0$ if $k \neq 0.$
\end{center}
\end{prop}
\noindent Instead of considering the second associated graded construction with respect to the maximal ideal $\m$ we could consider this construction with respect to every $\m$-primary ideal, but to simplify the exposition of the proofs we omit this generalization.

\begin{defi}
Let $(A,\m)$ be a Noetherian local ring, $I$ a proper ideal of $A$, $M$ a finitely generated $A$-module and $i$ an integer. Then we set
\begin{center}
$\locadeg{i}{I}{M}= \sum_{\p} \ell(H^{0}_{\p}(M_{\p})) \ee{i}{GG(A/\p)},$
\end{center}
where the sum is on the $i$-dimensional associated prime ideals of $M$.
\end{defi}

This definition is different from the one of \cite{MVY} (where the sum is taken over  the associated prime ideals of dimension $i + 1$, since they work with the projective dimension in a graded ring) but agrees with that of \cite{Vas}.  

\begin{rem}
We observe that the definitions of $M_{\leq i}$ and $M_{>i}$, which we gave in section 2.2 do not use that $M$ is graded.
By  Lemma \ref{LocCohom} and Proposition \ref{add} it is clear that 
$$\locadeg{i}{I}{M} = \gmult{i}{I}{M_{\leq i}}.$$
In the following we will use this characterization to compute the arithmetic degree. 
\end{rem}

\noindent This definition is also motivated by the following result, which is obtained by local duality as in the case where $I$ is $\m$-primary. 
 
\begin{prop}
Let $(A,\m)$ a d-dimensional Gorenstein local ring.
Then
\begin{center}

$\locadeg{i}{I}{M}= \ee{i}{GG(\ext{}{d-i}{M}{A})}.$
 
\end{center}

\end{prop}

Our definition agrees with the classical one given in the $\m$-primary case again applying Proposition \ref{clad}.
So, if $I=\m$ we obtain the usual notion of arithmetic degree, see Vasconcelos \cite{Vas}, p.~223.
We recall another result, which is straightforward from \cite{AMa}, Proposition 2.5, which we will use in the proof of the next theorem.

\begin{prop}\label{Sum}
Let $(A,\m)$ be a Noetherian local ring of dimension $d$, $I$ a proper ideal of $A$ and $M$ a \fg  $ A$-module. Then
\begin{center}
$\e{i}{\gr{I}{M}} = \displaystyle{\sum_{k=0}^{\dim M}} (\gmult{i}{I}{M})_{k}.$
\end{center} 
\end{prop}

Our next result is similar to a result of Sturmfels-Trung-Vogel \cite{STV} on arithmetic degree of an ideal and of its initial ideal. Since we are working with bigraded structures we have a more technical formulation, but the corollary clarifies this analogy.

\begin{teo} \label{Ineq}
Let $(A,\m)$ be a Noetherian local ring, $I$ a proper ideal in $A$, $M$ a finitely generated $A$-module and $\M$ the maximal homogeneous ideal in $\gr{I}{A}$.
Then

\begin{center}
$ \adeg{r}{\gr{\M}{\gr{I}{M}}} \geq \displaystyle{\sum_{k}} {(\locadeg{r}{I}{M})_{k}}. $
\end{center}

\end{teo}

Here, we consider $\gr{\M}{\gr{I}{M}}$ as a graded module on the graded ring $\gr{\M}{\gr{I}{A}}$, since it is already a quotient with respect to a homogeneous ideal of a polynomial ring on the field $A/\m$. 

In the proof we will need some results of Rees \cite{Rees} which we collect here for the convenience of the reader.

Let $(A,\m)$ be a Noetherian local ring and $\aid$ an ideal in $A$. We recall that the (extended) Rees algebra with respect to $\aid$ is the subring of $A[T, T^{-1}]$ given by $R(A,\aid) = \bigoplus_{r \in \mathbb{Z}}  \aid^{r} T^{r}$ where $\aid^{r} = A$ for $r \leq 0$. In the same way given a \fg $A$-module $M$ we set $R(M,\aid) = R(M) = \bigoplus_{r \in \mathbb{Z}} \aid^{r}T^{r}
M$.
Then, given an ideal $\bb \subseteq A$, we set $\bb^{*} = \bb A[T, T^{-1}] \cap R(A,\aid)$.
The following hold:
\begin{itemize}
\item[(i)] Set $A^{'}=A/\bb$ and $a^{'} = \aid + \bb / \bb \subseteq A^{'}$. Then $R(A^{'}, \aid^{'}) \simeq R(A,\aid) / \bb^{*} $.
\item[(ii)] If $\bb_{1}, \bb_{2}$ are two ideals in $A$, then $(\bb_{1} \cap \bb_{2})^{*} = \bb_{1}^{*} \cap \bb_{2}^{*} $.
\end{itemize}

From now on we shall use the symbol $T^{-1}$ to denote the map between $R(A,\aid)$-modules induced by the multiplication by $T^{-1}$.

\begin{proof}[Proof of Theorem \ref{Ineq}.]

We start with the following short exact sequence:
\begin{equation}\label{M<i}
0 \lgra M_{\leq r} \lgra M \lgra M_{>r} \lgra 0.
\end{equation}

\noindent Passing to the Rees algebra and setting $N_{0}=\bigoplus_{n} R(M)_{n} \cap M_{\leq r}$ gives the short exact sequence
\begin{center}
$ 0 \lgra N_{0} \lgra R(M) \lgra R(M_{>r}) \lgra 0.$
\end{center}

\noindent Since $(N_{0})_{n} = I^{n}M\cap M_{\leq r} $ and $R(M_{\leq r})_{n} = I^{n}M_{\leq r}$, we have $R(M_{\leq r}) \subseteq R(M)$. 
So $R(M_{\leq r})$ is a $R(A)$-submodule of $R(M)$ and we observe that $\dim R(M_{\leq r}) = \dim M_{\leq r} + 1 \leq r + 1 $, hence $R(M_{\leq r}) \subseteq R(M)_{\leq r + 1}$.

\noindent Let $P=R(M)_{\leq r + 1}/ R(M_{\leq r}) $ and $U, V$ be the kernel and the cokernel of the map  $T^{-1}\colon P(1) \lgra P$ respectively and set $G^{'} = R(M)_{\leq r+1} / T^{-1} R(M)_{\leq r+1}$.
\noindent Using the diagram

$$ \diagram
        &                                 &                                             &           0 \dto                     &    \\
        & 0 \dto                          &   0  \dto                           &      U  \dto        &   \\
0 \rto  & R(M_{\leq r})(1) \dto^{T^{-1}} \rto & R(M)_{\leq r+1} (1) \dto^{T^{-1}} \rto  & P(1) \dto^{T^{-1}} \rto & 0 \\
0 \rto  & R(M_{\leq r}) \rto \dto      & R(M)_{\leq r+1} \rto \dto                               & P             \rto \dto     & 0 \\
        & \gr{I}{(M_{\leq r})}\rto \dto & G^{'} \rto \dto & V \rto \dto  & 0 \\  
        &     0                    &     0       & 0       &   \\  
   \enddiagram
$$

\noindent the snake-lemma yields an exact sequence

\begin{equation}\label{G'}
0 \lgra U \lgra \gr{I}{(M_{\leq r})} \lgra G^{'} \lgra V \lgra 0 .
\end{equation}

{\bf Claim 1}: $G^{'} \subseteq (\gr{I}{M})_{\leq r} $.

\noindent We prove this in the same way as we have proved  $R(M_{\leq r}) \subseteq R(M)_{\leq r + 1}$.

\noindent At first we show  that $G^{'} \subseteq \gr{I}{M}$. Consider the exact sequence 

\begin{center}
$0 \lgra R(M)_{\leq r+1} \lgra R(M) \lgra R(M)_{> r} \lgra 0 $.
\end{center}

\noindent Multiplying by $T^{-1}$, we obtain a diagram similar to the previous one. But now $${\mathrm{Ker}}(R(M)_{> r+1}\stackrel{T^{-1}}{\lgra} R(M)_{> r + 1} )=0,$$ since $T^{-1}$ is a nonzerodivisor in $R(M)$. Hence, again by the snake-lemma, we have  $G^{'} \subseteq \gr{I}{M}$. Finally, to finish the proof of our claim, it remains to be shown that $\dim_{\gr{I}{A}} G^{'} \leq r $. Indeed, if $\varphi \colon R(A) \lgra R(A)/T^{-1}R(A) \simeq \gr{I}{A}$ we observe that 
$$\varphi ^{-1} (\annu{\gr{I}{A}}{G^{'}}) = \annu{R(A)}{G^{'}}.$$
Then 
\begin{center} $ R(A)/ \annu{R(A)}{G^{'}} \simeq \gr{I}{A}/ \annu{\gr{I}{A}}{G^{'}}$, \end{center}
\noindent so 
\begin{center}
$ \dim _{\gr{I}{A}} G^{'} = \dim _{R(A)} G^{'} = \dim R(M)_{\leq r +1 } -1 \leq r $
\end{center}

\noindent since $T^{-1}$ is a nonzerodivisor in $R(M)$.
 
\noindent We want to prove that $U$ and $V$ have the same cycle on $\mathrm{Spec}(\gr{I}{A})$. The idea is using, as in Theorem 1.2.6 of \cite{FOV}, the lemma of Artin-Rees and working on $P =R(M)_{\leq r+1} / R(M_{\leq r})$.

{\bf Claim 2}: $R(M)_{\leq r+1} \subseteq N_{0} $.

\noindent Let $x \in R(M)_{\leq r+1} \subseteq R(M)$. We can assume that $ x \neq 0 $ and that $x$ is a  homogeneous element of degree $n \in \mathbb{Z}$, hence $x \in I^{n}M$. Consider $\bb = \annu{A}{x}$. By the definition of $\bb^{*}$ we get $\bb ^{*} = \annu{R(A)}{x}$.
Set $A^{'} = A/\bb$, $I^{'} = I + \bb / \bb$. Then, since $x \in R(M)_{\leq r + 1}$, we have
\begin{center}
$\dim A/\bb + 1 = \dim R(A^{'},I^{'}) = \dim R(A,I)/\bb^{*} \leq r + 1$.
\end{center}
Hence $x \in (N_{0})_{n},$ therefore  $R(M)_{\leq r + 1} \subseteq N_{0} = \oplus_{n} (I^{n}M \cap M_{\leq r})$ and
\begin{center}
$ R(M_{\leq r}) \subseteq R(M)_{\leq r + 1} \subseteq N_{0} $
\end{center}and
\begin{center}
$P=  R(M)_{\leq r + 1} / R(M_{\leq r}) \subseteq N_{0} / R(M_{\leq r}) = \bigoplus_{n} (I^{n}M \cap M_{\leq r})T^{n} / I^{n}T^{n} M_{\leq r}$.
\end{center}

\noindent By the lemma of Artin-Rees the last quotient is annihilated by $T^{-k}$ for $k >> 0$, so the same holds for $P$ and this proves that the cycles of $V$ and $U$ are the same.

Using the equality of the cycles, the exact sequence (\ref{G'}) and Claim 1 we conclude that

\begin{center}
$\emult{}{\M}{\gr{I}{(M_{\leq r})}} = \emult{}{\M}{G^{'}} \leq \emult{}{\M}{(\gr{I}{M})_{\leq r}}.$
\end{center}

If we substitute in the above  argument the  sequence (\ref{M<i}) by the short exact sequence

\begin{center}
$ 0 \lgra (\gr{I}{M})_{\leq r} \lgra \gr{I}{M} \lgra (\gr{I}{M})_{>r} \lgra 0,$
\end{center}

we get
\begin{center}
$\e{}{\gr{\M}{(\gr{I}{M}_{\leq r})}} \leq \e{}{(\gr{\M}{\gr{I}{M}})_{\leq r}}.$
\end{center}

Putting together the last inequality and using Proposition \ref{Sum}, we obtain our thesis.

\end{proof}

\begin{cor}
Let $(A,\m)$ a Noetherian local ring, $I$ a proper ideal  and $M$ a \fg $ A$-module. Then
\begin{center}
$\adeg{i}{\gr{I}{M}} \geq \adeg{i}{M}.$
\end{center}
Here the arithmetic degree is considered with respect to the maximal (homogeneous) ideal. 
\end{cor}
\begin{proof}
Since $\adeg{i}{\gr{I}{M}} = \emult{}{\M}{(\gr{I}{M})_{\leq i}}$, by the proof of the Theorem \ref{Ineq} and by Proposition \ref{Sum}  we have
$$ \adeg{i}{\gr{I}{M}} = \emult{}{\M}{(\gr{I}{M})_{\leq i}} \geq$$
$$ \emult{}{\M}{\gr{I}{(M_{\leq i})}}  = \sum_{k=0}^{\dim M} (\locadeg{i}{I}{M})_{k}.$$

Using the same argument of Proposition 2.3 of \cite{AMa} we can show that the double sum transform of the Hilbert function of the bigraded module $GG(M)$ is
$$ H^{(1,1)}_{I,M}(i,j) =\ell (M/(\m^{i+1} + I^{j+1})M) + $$ 
$$ + \ell ( \bigoplus_{k=0}^{j} (I^{k}M \cap (\m^{i+1} + I^{k+1})M/(\m^{i+1}I^{k} + I^{k+1})M) $$

To prove  $\sum_{k}(\gmult{i}{I}{M_{\leq i}})_{k} \geq \sum_{k}(\gmult{i}{\m}{M_{\leq i}})_{k}$ it is sufficient to observe that $H^{(1,1)}_{I,M}(i,i) \geq H^{(1,1)}_{\m,M}(i,i)$ using the above description.

By Proposition \ref{clad} we can conclude.

\end{proof}
\begin{cor}
Let $(A,\m)$ be an equidimensional Noetherian local ring and $I$ a proper ideal.  If $A$ has embedded associated primes of dimension $i$, then $\gr{I}{A}$ has embedded associated primes of dimension $i$.\end{cor} 

\begin{proof}
It is well known that $\dim{A} = \dim \  {\gr{I}{A}}$ and that $\gr{I}{A}$ is equidimensional if $A$ is (see e.g. \cite{Fu} p.~436). So, since  $i \neq \dim{A}$ by the previous corollary
\begin{center}
   $\adeg{i}{\gr{I}{A}} \geq \adeg{i}{A} > 0,$ 
\end{center}
which implies that $\gr{I}{A}$ has embedded components.
\end{proof}

\end{section}

{\bf{Acknowledgment.}} I would like to thank R. Achilles for his warm encouragment and for having introduced me to this subject. Investigation partially supported by the MIUR and the University of Bologna, funds for Selected Research Topics.


\begin{thebibliography}{99}

\bibitem{AAl}
Achilles,~R.; Aliffi,~D.
{\it SEGRE: a script for the REDUCE package CALI}
{(1999-2001). Available at {\texttt{http://www.dm.unibo.it/\~{}achilles/segre.}}}

\bibitem{AcAv}
Achilles,~R.; Avramov,~L.
{\it Relations between properties of a ring and of its associated graded ring.}
{Seminar D. Eisenbud/ B. Singh/ W. Vogel, Vol. 2, pp. 5---29 Teubner-Text Math., \textbf{48}, Teubner, Leipzig 1982.}


\bibitem{AMa}
Achilles,~R.; Manaresi,~M.
{\it Multiplicities of a bigraded ring and intersection theory.}
{Math. Ann. \textbf{309} (1997), no.~4, 573--591.}

\bibitem{BM}
Bayer,~D.; Mumford,~D.
{\it What can be computed in algebraic geometry?}
{Computational Algebraic Geometry and Commutative Algebra (Cortona, 1991), 1--48, Sympos. Math. XXXIV, Cambridge Univ. Press, 1993.}


\bibitem{BroSh}
Brodmann,~M.P.; Sharp,~R.Y.
{\it Local cohomology: an algebraic introduction with geometric applications.}
{Cambridge Studies in Advanced Mathematics, \textbf{60}. Cambridge University Press, Cambridge, 1988.}

\bibitem{Cocoa}
Capani,~A.; Niesi,~G.; Robbiano,~L.
{\it CoCoA, a system for doing Computations in Commutative Algebra.}
{Available via anonymous ftp from: {\texttt{cocoa.dm.unige.it.}}}


\bibitem{Eis}
Eisenbud,~D.
{\it Commutative algebra. With a view toward algebraic geometry.}
{Graduate Texts in Mathematics, \textbf{150}. Springer-Verlag, New York, 1995.}

\bibitem{FOV}
Flenner,~H.; O'Carroll,~L.; Vogel,~W.
{\it Joins and intersection.}
{Springer Monographs in Mathematics. Springer-Verlag, Berlin, 1999.}

\bibitem{Fu}
Fulton,~W.
{\it Intersection theory.}
{Ergebnisse der Mathematik und ihrer Grenz\-ge\-bie\-te. 3. Folge. A Series of Modern surveys in Mathematics, 2}
{Springer-Verlag, Berlin, 1998.}

\bibitem{Red}
Gr\"abe,~H.-G.: 
{\it CALI - a REDUCE package for commutative algebra.}
{Version 2, (1995). Available through the REDUCE library {\texttt {redlib@rand.org}}.}

\bibitem{Mac}
Grayson,~D.~R., Stillman,~M.~E.
{\it Macaulay 2, a software system for research in algebraic geometry.}
{Available at {\texttt {http://www.math.uiuc.edu/Macaulay2/}}.}

\bibitem{Har}
Hartshorne,~R.
{\it Connectedness of the Hilbert scheme.}
{Inst. Hautes \'Etudes. Sci. Publ. Math. no.~29, 1966, 5--48.}

\bibitem{MVY}
Miyazaki,~C., Vogel,~W., Yanagawa,~K.
{\it Associated primes and arithmetic degrees.}
{J. Algebra \textbf{192} (1997), no.~1, 166--182.}

\bibitem{Rees}
Rees,~D.
{\it A note on form rings and ideals.}
{Mathematika \textbf{4} (1957), 51--60.}

\bibitem{Ru}
Ruppert,~M.~W.
{\it Embedded primes and the associated graded ring.}
{Arch. Math. (Basel)  \textbf{65}, (1995), no.~6, 501--504.}


\bibitem{STV}
Sturmfels,~B., Trung,~N.~V., Vogel,~W.
{\it Bounds on degrees of projective schemes.}
{Math. Ann. \textbf{302}, (1995), no.~3, 417--432.}

\bibitem{Trung1}
Trung,~N.~V.
{\it Positivity of mixed multiplicities.}
{Math. Ann. \textbf{319}  (2001), no.~1, 33--63.}


\bibitem{VdW}
Van der Waerden,~B.~L.:
{\it On Hilbert's function, series of composition of ideals and a generalization of the theorem of Bezout.}
{Proc. Roy. Acad. Amsterdam \textbf{31} (1929), 749--770.}


\bibitem{Vas}
Vasconcelos,~W.~V.:
{\it Cohomological degrees of graded modules.}
{Six Lectures on commutative algebra. (Bellaterra, 1996), 345--392, Prog. Math., \textbf{166}, Birkh\"auser, Basel, 1998.}




\end{thebibliography}
\end{document}